\documentclass{tran-l}
 \newtheorem{thm}{Theorem}[section]
 
 \newtheorem{lem}[thm]{Lemma}
 
 \theoremstyle{definition}
 
 \theoremstyle{remark}
 
 \numberwithin{equation}{section}

\def\N{{\mathbb N}}

\def\P{{\mathbb P}}

\begin{document}

\title[Counting Primes in the Interval
$\big(n^2,(n+1)^2\big)$]
 {Counting Primes in the Interval
$\big(n^2,(n+1)^2\big)$}

\author{Mehdi Hassani}

\address{Institute for Advanced\\
Studies in Basic Sciences\\
P.O. Box 45195-1159\\
Zanjan, Iran.}

\email{mmhassany@srttu.edu}

\thanks{}



\date{}

\dedicatory{}

\commby{}
%%% ----------------------------------------------------------------------
\begin{abstract}
In this note, we show that there are many infinity positive
integer values of $n$ in which, the following inequality holds
$$
\left\lfloor\frac{1}{2}\left(\frac{(n+1)^2}{\log(n+1)}-\frac{n^2}{\log
n}\right)-\frac{\log^2 n}{\log\log
n}\right\rfloor\leq\pi\big((n+1)^2\big)-\pi(n^2).
$$
\end{abstract}

\maketitle

%\tableofcontents
\section{Introduction}
Considering Euclid's proof for the existence many infinity primes,
we can get the following inequality for many infinite values of
$n$:
$$
1\leq\pi\big((n+1)^2\big)-\pi(n^2),
$$
in which $\pi(x)=\#[2,x]\cap\P$, and $\P$ is set of all primes.
Now, we have some strong results, which allow us to change 1 in
left hand side of above inequality by a nontrivial one. In fact,
we show that there are many infinity positive integer values of
$n$, in which the following inequality holds:
$$
\left\lfloor\frac{1}{2}\left(\frac{(n+1)^2}{\log(n+1)}-\frac{n^2}{\log
n}\right)-\frac{\log^2 n}{\log\log
n}\right\rfloor\leq\pi\big((n+1)^2\big)-\pi(n^2).
$$
This is the result of an unsuccessful challenge, for proving the
old-famous conjecture, which asserts for every $n\in\N$, the
interval $\big(n^2,(n+1)^2\big)$ contains at least a prime.
Surely, Prime Number Theorem \cite{davenport}, suggests a few more
number of primes as follows:
$$
F(n)\sim\frac{1}{2}\left(\frac{(n+1)^2}{\log(n+1)}-\frac{n^2}{\log
n}\right)\hspace{10mm}(n\rightarrow\infty),
$$
in which $F(n)$ is the number of primes in
$\big(n^2,(n+1)^2\big)$. This asymptotic relation, led us to make
some conjectures on the bounding $F(n)$.\\\\
\textbf{Conjecture 1.} For every $n\geq 5$, we have
$$
F(n)<\frac{1}{2}\left(\frac{(n+1)^2}{\log(n+1)}-\frac{n^2}{\log
n}\right)+\log^2 n \log\log n.
$$
This conjecture has been checked by Maple for all $5\leq n\leq
10000$.\\\\
\textbf{Conjecture 2.} For every $n\geq 3$, we have
$$
\frac{1}{2}\left(\frac{(n+1)^2}{\log(n+1)}-\frac{n^2}{\log
n}\right)-\frac{\log^2 n}{\log\log n}-1<F(n).
$$
This conjecture has been checked by Maple for all $3\leq n\leq
10000$. Also, as mentioned above, we show that for many infinity
positive integer values of $n$, the truth of this conjecture
holds. To do this, we need the following sharp bounds for the
function $\pi(x)$ (see \cite{dusart}):
\begin{equation}\label{L}
L(x)=\frac{x}{\log x}\Big(1+\frac{1}{\log
x}+\frac{1.8}{\log^2x}\Big)\leq\pi(x)\hspace{10mm}(x\geq 32299),
\end{equation}
and
\begin{equation}\label{U} \pi(x)\leq U(x)=\frac{x}{\log
x}\Big(1+\frac{1}{\log
x}+\frac{2.51}{\log^2x}\Big)\hspace{10mm}(x\geq 355991).
\end{equation}

\section{Main Result}
\begin{lem} For every $n\geq 2$, we have
$$
\frac{n^2}{2\log n}+4-\frac{9}{\log
9}-\sum_{k=3}^{n-1}\frac{\log^2 k}{\log\log k}<\frac{n^2}{2\log
n}\left(1+\frac{1}{2\log n}+\frac{9}{20\log^2 n}\right).
$$
\end{lem}
\begin{proof} For every $n\geq 2$, consider the following
inequality
$$
\frac{n^2}{4\log^2 n}+\frac{9n^2}{40\log^3
n}+\sum_{k=3}^{n-1}\frac{\log^2 k}{\log\log k}>4-\frac{9}{\log 9}.
$$
Note that the left member of it, is positive and the right member
is negative. So, clearly it holds for every $n\geq 2$.
\end{proof}
\begin{lem}\label{pi-n2} For every $n\geq 180$, we have
$$
\frac{n^2}{2\log n}+4-\frac{9}{\log
9}-\sum_{k=3}^{n-1}\frac{\log^2 k}{\log\log k}<\pi(n^2).
$$
\end{lem}
\begin{proof} Putting $x=n^2$ in (\ref{L}), for $n\geq
180=\left\lceil\sqrt{32299}\right\rceil$ we obtain
$$
\frac{n^2}{2\log n}\left(1+\frac{1}{2\log n}+\frac{9}{20\log^2
n}\right)<\pi(n^2).
$$
Considering this, with previous lemma, completes the proof.
\end{proof}
\begin{thm} For many infinity positive integer values of $n$,
the following inequality holds
$$
\left\lfloor\frac{1}{2}\left(\frac{(n+1)^2}{\log(n+1)}-\frac{n^2}{\log
n}\right)-\frac{\log^2 n}{\log\log
n}\right\rfloor\leq\pi\big((n+1)^2\big)-\pi(n^2).
$$
\end{thm}
\begin{proof} Reform the truth of lemma \ref{pi-n2}, as follows:
$$
\frac{1}{2}\left(\frac{n^2}{\log n}-\frac{3^2}{\log
3}\right)-\sum_{k=3}^{n-1}\frac{\log^2 k}{\log\log
k}<\pi(n^2)-\pi(3^2).
$$
This inequality yields the following one:
$$
\sum_{k=3}^{n-1}\left\lfloor\frac{1}{2}\left(\frac{(k+1)^2}{\log(k+1)}-\frac{k^2}{\log
k}\right)-\frac{\log^2 k}{\log\log
k}\right\rfloor<\sum_{k=3}^{n-1}\pi\big((k+1)^2\big)-\pi(k^2),
$$
which holds for all $n\geq 180$. Now, we note that terms under
summations, in both sides are non-negative integers and this
completes the proof\footnote{In fact we can show that if $a_n$ and
$b_n$ are two non-negative integer sequences, with $\sum_{n=n_0}^N
a_n<\sum_{n=n_0}^N b_n$, then we have $\#\{n|~a_n\leq
b_n\}=\aleph_0$.}.
\end{proof}
However, this challenge was unsuccessful for proving the relation
$$
\{n|~\big(n^2,(n+1)^2\big)\cap\P\neq\phi\}=\N,
$$
but it seems that it can be useful for improving it. To see this,
let
$$
g(n)=\#\big\{t~|~t\in\N,~t\leq n,~\P\cap
\big(t^2,(t+1)^2\big)\neq\phi\big\}.
$$
Clearly, $\lim_{n\rightarrow\infty}g(n)=\infty$ and $g(n)\leq n$.
Note that $g(n)=n$ is above mentioned open problem. A lower bound
for $g(n)$ is the following bound, which can yield by considering
previous theorem for every $n\geq 597$;
$$
g(n)\geq M(n),
$$
in which
$$
M(n)=\max_{m}\left\{\sum_{k=597}^{n}\left\lfloor\frac{1}{2}\left(\frac{(k+1)^2}{\log(k+1)}-\frac{k^2}{\log
k}\right)-\frac{\log^2 k}{\log\log
k}\right\rfloor\leq\sum_{k=m}^{n}U\big((k+1)^2\big)-L(k^2)\right\}.
$$
Clearly, if $n\rightarrow\infty$, then we have
$$
M(n)=O(n).
$$
Also, we have the following conjecture on the size of $M(n)$:\\\\
\textbf{Conjecture 3.} For every $\epsilon>0$ there exists
$n_{\epsilon}\in\N$ such that for all $n>n_{\epsilon}$ we have
$$
M(n)>(1-\epsilon)n.
$$

% ------------------------------------------------------------------------

\end{document}